\documentclass[12pt,oneside]{amsart}
\usepackage{latexsym}
\usepackage{amsmath,amssymb,amsthm}
\usepackage{amscd}
\usepackage[dvips]{graphics}
\newcommand{\Z}{\mathbb{Z}}
\newcommand{\R}{\mathbb{R}}
\newcommand{\C}{\mathbb{C}}

\theoremstyle{remark}

\newtheorem*{ack}{Acknowledgements}

\begin{document}
\title{Seidel's Mirror Map for the Torus}
\author{Eric Zaslow}
\begin{abstract} 
Using only the Fukaya category and the monodromy around large complex
structure, we reconstruct the mirror map in the case of a symplectic torus. 
This realizes an idea described by Paul Seidel.
\end{abstract}
\maketitle

\section{Introduction}
\label{intro}

Paul Seidel had the following idea for
recovering the mirror map purely from the Fukaya category.\footnote{The
idea described was told in a private communication to the author. 
This may have been implicit in the
works of Fukaya and/or in the minds of others in the field.}  Start with a symplectic
Calabi-Yau $X$ and its family of complex structures, and
assume it has a projective mirror manifold $Y$ with a family of symplectic
structures, and that
Kontsevich's conjecture holds:  $DFuk(X) \cong D(Y),$
where $DFuk(X)$ is the Fukaya category  of $X$ (i.e. the bounded derived category
constructed from the Fukaya $A_\infty$ category)
and $D(Y)$ is the bounded derived category of coherent
sheaves on $Y$.  Then the
homogeneous coordinate ring
%of holomorphic functions
on $Y$ is given by
$$\mathcal R = \bigoplus_{k=0}^\infty \Gamma(\mathcal O_Y(k)) =
\bigoplus_{k=0}^\infty {\rm Hom}_{DFuk(X)}(\psi(\mathcal O), \psi(\mathcal O(k))),$$
where $\psi$ is the equivalence of categories.
The term on the right can be evaluated solely in $DFuk(X)$, and
thus the complex projective variety $Y$ can be recovered.  The dependence
of this construction on the symplectic structure of $X$ defines the mirror
map.\footnote{The case of Fano varieties is being considered in \cite{abd}.}

Let $S \equiv \psi(\mathcal O)$ be the object dual to the structure sheaf of $Y,$
conjecturally the Lagrangian section of the Lagrangian torus fibration (cf. \cite{syz}).
We will often
equate a geometric Lagrangian submanifold with the object in $DFuk$ which
it defines, including, if necessary, additional data such as grading and
local system.
Recall \cite{horjakontseidthom} that on the complex
structure moduli space of $X,$ monodromies act by symplectomorphisms,
which define autoequivalences of $DFuk(X)$ (we use the same notation
for a symplectomorphism and the autoequivalence it induces) and that
the monodromy $\rho$ around the large complex
structure limit point is mirror to the autoequivalence of $D(Y)$ defined by
$\mathcal E \rightarrow \mathcal E\otimes \mathcal O(1).$
We define $L_k$ by $L_k \equiv \rho^k S.$  Note $S = L_0$
and $L \equiv L_1$ is dual to $\mathcal O(1).$  In fact, $L_k = \psi(\mathcal O(k)),$ 
so we wish to compute $\bigoplus_i {\rm Hom}_{DFuk(X)}(S,L_k).$
In order to interpret this as a ring, we must identify ${\rm Hom}(L_k,L_{k+l}))$
with ${\rm Hom}(S,L_l)$
(we hereafter drop the $DFuk(X)$ subscript),
and to do so we use the symplectomorphism $\rho^{-k}.$

In this note we will compute $\mathcal R$ in the case where $X$ is a
symplectic two-torus and derive the mirror map.\footnote{The result is guaranteed to
be correct here, since Kontsevich's conjecture has been proven
in this example \cite{pz}.}
Without knowing the mirror map, we can still say
that $Y$ is some elliptic curve, and thus
has a projective embedding as a cubic curve.  Then $\mathcal O(1)$
is a line bundle of degree three on $Y,$ so its mirror must have intersection
three with $S.$  Taking the base section $S$ to be the $x$-axis in the
universal cover $\R^2,$ we have that $L$ is a
line of slope three.  So we put $\rho = \gamma^3,$ where $\gamma$ is
a minimal Dehn twist, and note that $\rho$ is maximally unipotent.
For simplicity we take $S$ (and therefore $L$) to have trivial local systems
and to pass through lattice vectors, but our results do not depend on this choice.
The data of $S$ and $\rho$ now allows us to calculate $\mathcal{R}.$

\begin{ack}
I would like to thank Paul Seidel for communicating his ideas freely.  Thanks to
The Fields Institute for hosting me during this project.  This work was
supported in part by a Clay Senior Scholars fellowship and by NSF grant
DMS--0405859.
\end{ack}

\section{Computation}
\label{computation}

We define $X = \R^2/\Z^2$ with $\omega = \tau dx\wedge dy,$
$\tau\in\C,$ $Im(\tau)>0.$  The category constructed from Fukaya's
$A_\infty$ category in this case was described explicitly in \cite{pol1,pol2,pz},
and we refer the reader to those papers for details.
As discussed above, we have
$L_k = \{ (t,3kt) \mod \Z^2 : t \in \R \},$ and we define
its grading $\alpha = \tan^{-1}(k) \in [0,\pi/2).$  We define
$X_i = (i/3,0) \in {\rm Hom}(S,L),$
$Y_i = (i/6,0) \in {\rm Hom}(S,L_2),$
and $Z_i = (i/9,0) \in {\rm Hom}(S,L_3),$ where $i$ is taken mod $3,$ $6,$ and
$9,$ respectively.  In the sequel,
when we write an equation like $X_1X_2 = ...,$ the $X_2$
is understood to live in ${\rm Hom}(L_1,L_2)$ through $\rho.$
Explicitly, $\rho(x,y) = (x,y+3x);$ indeed $\rho^*\omega = \omega.$

Let us compute the products $X_iX_j.$  The Fukaya
category for this example was discussed in \cite{kont, pz}.
The basic computation is $X_0X_1.$  The minimal triangle (holomorphic map)
appearing in the product
connects the points $X_0 = (0,0),$ $\rho(X_1) = X_1 = (1/3,1),$ and $Y_1 = (1/6,0)$
and has symplectic area $(1/2)(1/6)(1)\tau.$  Multiples and translates of
this triangle are relevant to other products.  Multiples by $6n$ have the
same endpoints and contribute to the same product, with
area $(1/2)(n+1/6)(6n+1)\tau.$  The coefficient of $Y_1$ in
$X_0X_1$ is thus $A_1
\equiv \sum_n \exp [ i \pi 6 \tau (n+1/6)^2 ] = \theta[1/6,0](6\tau,0).$\footnote{We
recall the definition $\theta[a,b](\tau,z)=
\sum_{n\in\Z}\exp[i\pi\tau(n+a)^2+2\pi i (n+a)(z+b)].$}
Defining $A_k := \theta[k/6,0](6\tau,0),  k \in \Z/6\Z,$ and noting $A_k = A_{6-k},$
we get the following relations:
\begin{equation}
\label{aprod}
X_iX_j = \sum_{k=0}^1 A_{i-j+3k}Y_{i+j+3k}.
\end{equation}
The right hand side of this equation makes sense with $i,$ $j$ defined mod $3$.
Commutativity is easily shown to follow from the relations among the $A_k.$

Next we compute $Y_iX_j.$  Starting with $Y_1 X_1,$ the
minimal triangle has vertices $Y_1=(1/6,0),$ $\rho^2(X_1)=X_1=(1/3,1),$ and
$Z_2 = (2/9,0),$ with area $(1/2)(1/18)(1)\tau.$  Odd multiples 
(with left endpoint fixed) 
and translates of this triangle are relevant to $Y_iX_j$ with $i$ odd;
even multiples and translates to $i$ even.  Multiples by $18n$ have
the same endpoints.  Therefore
$Y_1X_1 = B_1 Z_2 + B_7 Z_5 + B_{13}Z_8,$
where $B_k = \sum_n \exp [ i \pi 18 \tau (n+k/18)^2 ] = \theta[k/18,0](18\tau,0).$
Note $B_k = B_{18-k}$ and $k$ is defined mod $18.$
As an example of another product, the third multiple of the minimal
triangle has endpoints $Y_1=(1/6,0),$ $X_2=(2/3,3),$ $Z_3=(1/3,0),$
thus $Y_1X_2 = B_3Z_3 +....$
Collecting results, we find
\begin{equation}
\label{bprod}
Y_i X_j = \sum_{k=0}^2 B_{2j-i+6k}Z_{i+j+3k}.
\end{equation}

\section{Commutativity and Associativity}
\label{assoc}

Associativity in the (derived or cohomological) Fukaya category
follows from general grounds, and in the case of the torus amounts
to an equality obtained from expressing the area of a non-convex quadrangle
by splitting it into triangles in two different ways.  (This was noted, for
example, in Section 2 of
\cite{pol1}.)  It also amounts to relations among the $A_k$ and $B_k,$
which we describe presently.

As for commutativity, this follows from the existence of a robust family of
anti-symplectomorphisms.  For example, in considering the products $X_0Y_k,$
one must count (among other things) triangles with vertices $X_0,$
$\rho(Y_k),$ and $Z_{k}$ 
arranged in clockwise orientation and with sides of appropriate slope.
Now consider the map $\varphi:$
$$(x,y)\mapsto \left(\frac{1}{2}x - \frac{7}{18}y + \frac{1}{9}k\;, -2y\right).$$
We note $\varphi(X_0) = Z_k,$ $\varphi(\rho(Y_k)) =
X_0 = \rho^2(X_0),$ and $\varphi(Z_k) = Y_k.$   Further, since $\varphi$ is an anti-symplectomorphism,
i.e. $\varphi^*\omega = -\omega,$ it preserves areas
and reverses the orientation and thus changes
the order in which the vertices appear on the outside of the triangle.  Thus
$Y_k, \rho^2(X_0), Z_k$ are oriented clockwise in the image triangle,
which has the same area as the original.  This proves commutativity among
products $X_0 Y_k.$  Translations of $\varphi$ suffice for proving
commutativity for $X_j Y_k.$  Products $X_iX_j$ were already seen to be
commutative, and this is all that we will require for our purposes. 
In short, commutativity follows from anti-symplectomorphisms mapping
vertices $(X,\rho^n Y,Z)$ to $(Z,\rho^m X,Y)$ in holomorphic triangles.
It is not clear (to the author)
why commutativity should hold in a general symplectic
manifold.

We now return to an explicit description of the associativity constraint.
We will make use of the following identity, which
follows from the addition formula II.6.4 of \cite{mumford}:
\begin{align}
\label{mum}
&\theta\left[\frac{a}{n},0\right](n\tau,0)\; \theta\left[\frac{b}{nk},0\right](nk\tau,0)=\\
\nonumber
&
\sum_{\epsilon=0}^{k}
\theta\left[\frac{b-ka+kn\epsilon}{k(k+1)n},0\right](k(k+1)n\tau,0)\;
\theta\left[\frac{a+b+kn\epsilon}{(k+1)n},0\right]((k+1)n\tau,0).
\end{align}
When $n=6$ and $k=3$ this gives us formulas for $A_aB_b.$
Defining $C_c = \theta[c/24,0](24\tau)$ and $D_d = \theta[d/72,0](72\tau),$
we have
\begin{equation}
\label{simple}
A_aB_b = \sum_{\epsilon=0}^{3}C_{a+b+18\epsilon}D_{b-3a+18\epsilon}.
%=\overline{C_{a+b}D_{b-3a}},
\end{equation}
%where the overline indicates the sum over elements of the orbit
%of the subscripts $(a+b,b-3a)$ in $\Z_{24}\times\Z_{72}$ under translations
%by $(18,18).$\footnote{Be careful:  although $C_c = C_{24-c},$ it is not
%necessarily true that $\overline{C_cD_d} = \overline{C_{24-c}D_d}.$}
This formula suffices for proving some of the equivalences necessary for
showing associativity.  Others follow from further application of (\ref{mum}).
%For others, further application of (\ref{mum})
%or other tricks are required.  We note $C_cD_d = \overline{E_{c+d}F_{d-3c}},$
%where $E$ and $F$ are defined analogously to $C$ and $D,$ only
%for $96\tau$ and $288\tau,$ respectively, and the overline indicates
%the orbit under $(72,72)\in \Z_{96}\times\Z_{288}.$

For example, one wants to show that $(X_0^2)X_1 = X_0(X_0X_1).$
This amounts to $(A_0Y_0+A_3Y_3)X_1=
X_0(A_1Y_1+A_2Y_4).$  Using commutativity and the products (\ref{aprod}),
then equating coefficients on $Z_k,$ gives the conditions
\begin{eqnarray*}A_0B_2+A_3B_7&=&A_1B_1+A_2B_8,\\
A_0B_8+A_3B_1&=&A_1B_5+A_2B_4,\\
A_0B_4+A_3B_5&=&A_1B_7+A_2B_2.
\end{eqnarray*}
The first and third relations follow immediately from (\ref{simple}).
The second equation is most easily seen by rewriting the right hand
side as $A_{-1}B_5 + A_{-2}B_{-4}.$
Proceeding in this manner, one can prove well-definedness of $X_iX_jX_k.$

Again, associativity follows from quadrilateral dissection, or on general
grounds for the Fukaya category, and our philosophy here should be
to think of these identities as following from the associativity constraints.
In either case, we will use the explicit expressions derived here.

\section{Relations}
\label{rels}

One finds that the number of degree two polynomials in the three variables $X_i$
equals exactly the number of $Y_k,$ and in fact since $A_0A_1-A_2A_3 \neq 0$
one finds that
the $Y_k$ can be written in terms of products $X_iX_j,$ and vice versa, so
there are no relations in $\mathcal R$ at this degree.
At the next level, we have ten independent polynomials and nine $Z_k,$
so we expect a single relation.  Let us search for this relation.

Let $$\{X_0^3,X_1^3,X_2^3,X_0^2X_1,X_1^2X_2,X_2^2X_0,X_0^2X_2,
X_1^2X_0,X_2^2X_1,X_0X_1X_2\}$$ be a basis, with $e^I$ the $I$-th entry,
$I = 0...9.$  Using the product, we can write $e^I = \sum_k M_k{}^I Z_k.$
A relation $a$ has the form $\sum_I a_I e^I = 0,$ or
$\sum_k \left( \sum_I (M_k{}^I a_I)\right)  Z_k = 0.$  Since the $Z_k$ are
linearly independent generators of ${\rm Hom}(S,L_3)$ we have,
in matrix form $M\cdot a = 0,$ or $a\in {\rm Ker}(M).$  $M$ is a $9\times 10$
matrix, so the kernel should be one-dimensional, and we can take
$a_I = c (-1)^I {\rm det}(M_I),$ where $M_I$ is $M$ with the $I$-th column
removed and $c\neq 0$ is any constant.

Using the products found in Section \ref{computation}, one finds
$$M = \left(\begin{array}{cccccccccc}
p&q&q&0&0&0&0&0&0&u\\
0&0&0&r&t&s&0&0&0&0\\
0&0&0&0&0&0&t&r&s&0\\
q&p&q&0&0&0&0&0&0&v\\
0&0&0&s&r&t&0&0&0&0\\
0&0&0&0&0&0&s&t&r&0\\
q&q&p&0&0&0&0&0&0&v\\
0&0&0&t&s&r&0&0&0&0\\
0&0&0&0&0&0&r&s&t&0
\end{array}\right),$$
where
$$\begin{array}{lll}
p=A_0B_0+A_3B_9\qquad&r=A_0B_2+A_3B_7\qquad&u=A_2B_0+A_1B_9\\
q=A_0B_6+A_3B_3\qquad&s=A_0B_8+A_3B_1\qquad&v=A_2B_6+A_1B_3\\
{}&t=A_0B_4+A_3B_5.&{}
\end{array}$$
Up to a common multiple, one finds
$a \sim ( (p+q)u - 2qv, pv - qu, pv-qu, 0,0,0,0,0,0,2q^2 - pq - p^2 ).$
In fact, $u=v$, which follows from associativity, or equivalently
the relations (\ref{mum}), so we can remove the common (nonzero)
factor of $p-q$ and take
$$a = (u,u,u,0,0,0,0,0,0,-2q-p).$$

If there are no other relations in the ring $\mathcal R,$ then
this single relation defines a cubic curve in the Hesse family as
$$a_0 X_0^3 + a_1 X_1^3 + a_2 X_2^3 + a_9 X_0X_1X_2 = 0.$$
The modular invariant is easily calculated in terms of
$z = -(1/3)a_9(a_0a_1a_2)^{-1/3}=\frac{2q+p}{3u}.$
Explicitly,
\begin{equation}
\label{jeq}
j(\tau) = -27z^3(z^3+8)^3(1-z^3)^{-3}.
\end{equation}
This equation, which should define the $j$-function of the mirror curve,
is written in terms of the symplectic parameter $\tau$ on the torus.
It therefore defines the mirror map, which in this example is known
to send the symplectic parameter $\tau$ to the modular parameter $\tau$
in the upper halfplane. 
So (\ref{jeq}) amounts to an identity in terms of the
variable $\tau,$ or more conveniently for us, $x = e^{i\pi \tau/18},$
and it remains to verify this relation.\footnote{We
ignore the possibility of further relations in $\mathcal R.$
This assumption is justified using the mirror equivalence,
but would be difficult to show working purely from the Fukaya side.}

The following identities follow directly from the definitions:
\begin{eqnarray*}
A_k &=& x^{3k^2} + \sum_{n=1}^\infty x^{3(6n+k)^2} + x^{3(6n-k)^2},\\
B_k &=& x^{k^2} +  \sum_{n=1}^\infty x^{(18n+k)^2} + x^{(18n-k)^2}.
\end{eqnarray*}
Recall that the $j$-invariant has the expansion
$$
j(x) = x^{-36} + 744 + 196884x^{36} + 21493760x^{72} + 864299970x^{108}
%+ 20245856256x^{144}
+ ....$$
These coefficients and more
can be corroborated order by order in the series expansion of the
right hand side of (\ref{jeq}).  
A more general proof may be found in \cite{bl}.
Of course, this had to be true, by the equivalence of categories already proven
in \cite{pz},
but our intent was to find this result working only from the Fukaya
category.\footnote{Perhaps one could invert this philosophy
and derive information about the Fukaya category from the known
mirror maps, in cases where computing products is formidable.} 
We find the computation a pleasant realization of Seidel's idea.

\bibliographystyle{abbrv}
%\bibliography{thesis}

\vskip 0.1in

%{\scriptsize John Loftin, Department of Mathematics, Columbia
%University, New York, NY  10027.
%(LOFTIN@MATH.COLUMBIA.EDU)}

%{\scriptsize Shing-Tung Yau, Department of Mathematics, Harvard
%University, Cambridge, MA  01238. (YAU@MATH.HARVARD.EDU) }

\noindent
{\scriptsize Eric Zaslow, Department of Mathematics,}
{\scriptsize Northwestern University,}
{\scriptsize Evanston, IL  60208.}\\
{\scriptsize (zaslow@math.northwestern.edu)}

\end{document}